\documentclass{amsart}

\usepackage{amssymb,amsthm,mathtools,bm,color,scalerel}
\usepackage{hyperref}

\usepackage[T1]{fontenc}

\usepackage{tikz}
\usetikzlibrary{positioning}
\usetikzlibrary{svg.path}
\definecolor{orcid_color}{HTML}{A6CE39}

\hypersetup{pdfborder={0 0 0}} 
\DeclareRobustCommand{\orcidicon}{%
	\raisebox{.2mm}{\scalerel*{%
	\begin{tikzpicture}[xscale=1,yscale=-1,transform shape]
	\filldraw[color=orcid_color] svg {M256,128c0,70.7-57.3,128-128,128C57.3,256,0,198.7,0,128C0,57.3,57.3,0,128,0C198.7,0,256,57.3,256,128z};
	\filldraw[color=white] svg {M86.3,186.2H70.9V79.1h15.4v48.4V186.2z} svg {M108.9,79.1h41.6c39.6,0,57,28.3,57,53.6c0,27.5-21.5,53.6-56.8,53.6h-41.8V79.1z M124.3,172.4h24.5
		c34.9,0,42.9-26.5,42.9-39.7c0-21.5-13.7-39.7-43.7-39.7h-23.7V172.4z} svg {M88.7,56.8c0,5.5-4.5,10.1-10.1,10.1c-5.6,0-10.1-4.6-10.1-10.1c0-5.6,4.5-10.1,10.1-10.1
		C84.2,46.7,88.7,51.3,88.7,56.8z};
	\end{tikzpicture}}{|}}%
}
\newcommand{\orcid}[1]{\href{https://orcid.org/#1}{\orcidicon}}

\newcommand{\arxiv}[1]{arXiv:\href{https://doi.org/10.48550/arXiv.#1}{#1}}

\newtheorem{theorem}{Theorem}
\newtheorem{lemma}[theorem]{Lemma}

\begin{document}

\title{Characterizing graphs with fully positive semidefinite $Q$-matrices}
\author[Hajime Tanaka]{Hajime Tanaka\,\orcid{0000-0002-5958-0375}}
\address{\href{https://www.math.is.tohoku.ac.jp/index.html}{Research Center for Pure and Applied Mathematics}, Graduate School of Information Sciences, Tohoku University, Sendai 980-8579, Japan}
\email{htanaka@tohoku.ac.jp}
\urladdr{https://hajimetanaka.org/}
\begin{abstract}
For $q\in\mathbb{R}$, the $Q$-\emph{matrix} $Q=Q_q$ of a connected simple graph $G=(V,E)$ is $Q_q=(q^{\partial(x,y)})_{x,y\in V}$, where $\partial$ denotes the path-length distance.
Describing the set $\pi(G)$ consisting of those $q\in \mathbb{R}$ for which $Q_q$ is positive semidefinite is fundamental in asymptotic spectral analysis of graphs from the viewpoint of quantum probability theory.
Assume that $G$ has at least two vertices.
Then $\pi(G)$ is easily seen to be a nonempty closed subset of the interval $[-1,1]$.
In this note, we show that $\pi(G)=[-1,1]$ if and only if $G$ is isometrically embeddable into a hypercube (infinite-dimensional if $G$ is infinite) if and only if $G$ is bipartite and does not possess certain five-vertex configurations, an example of which is an induced $K_{2,3}$.
\end{abstract}

\maketitle

\hypersetup{pdfborder={0 0 1}} 

\section{Introduction}

For $q\in\mathbb{R}$, the $Q$-\emph{matrix} $Q=Q_q$ of a connected simple graph $G=(V,E)$ with vertex set $V$, edge set $E$, and path-length distance $\partial$, is defined by
\begin{equation*}
	(Q_q)_{x,y}=q^{\partial(x,y)} \quad (x,y\in V),
\end{equation*}
where we set $0^0:=1$.
The $Q$-matrix arises in defining the $q$-\emph{deformed vacuum functional} (or the \emph{Gibbs functional}) on the adjacency algebra of $G$ and plays an important role in asymptotic spectral analysis of graphs from the viewpoint of quantum probability theory; see, e.g., \cite{HO2007B,HO2008TAMS}.
For the $q$-deformed vacuum functional to be a state, $Q_q$ must be positive semidefinite in the sense that
\begin{equation*}
	\sum_{x,y\in V} \!\overline{f(x)}f(y)\,q^{\partial(x,y)} \geqslant 0
\end{equation*}
for every $f\in C_0(V)$, the space of $\mathbb{C}$-valued functions on $V$ with finite supports;
in other words, the function $(x,y)\mapsto q^{\partial(x,y)}$ $(x,y\in V)$ is a positive definite kernel on $V$.
Hence it is fundamental to consider the following subset of $\mathbb{R}$:
\begin{equation*}
	\pi(G)=\{q\in\mathbb{R}:Q_q\succcurlyeq 0\},
\end{equation*}
where $\succcurlyeq 0$ means positive semidefinite.
The set $\pi(G)$ is nonempty since $0,1\in\pi(G)$ and is closed.
We have $\pi(G)=\mathbb{R}$ if $G$ is a singleton, so assume that $G$ has at least two vertices.
Then
\begin{equation*}
	\begin{bmatrix} 1 & q \\ q & 1 \end{bmatrix}
\end{equation*}
is a principal submatrix of $Q_q$ and must be positive semidefinite whenever $q\in \pi(G)$.
It follows that
\begin{equation*}
	\pi(G)\subset [-1,1].
\end{equation*}
In general, it is a difficult problem to describe $\pi(G)$.
At an \href{https://www.impan.pl/en/activities/banach-center/conferences/22-19thworkshop}{international workshop} held in B\k{e}dlewo, Poland, in August 2022, which the author attended virtually, Marek Bo\.{z}ejko asked if there is a characterization of graphs $G$ for which $\pi(G)=[-1,1]$.
This note aims to provide such a characterization, which is combinatorial in nature.
See also \cite{KOT2021SIGMA,Obata2007SM,Obata2011CM,Voit2021pre} and the references therein for more results on the $Q$-matrix.
When we consider the $q$-deformed vacuum functional, we must assume that the graph $G$ is locally finite, for otherwise the adjacency algebra is not defined in the first place.
However, in this note, we will not make this assumption and will focus only on the positive semidefiniteness of the $Q$-matrix.

In Section \ref{sec: proofs}, we first prove the following lemma, which is of independent interest:
\begin{lemma}\label{-1 is in pi(G)}
For a connected simple graph $G$, the following are equivalent:
\begin{enumerate}
\item $-1\in \pi(G)$.
\item $G$ is bipartite.
\end{enumerate}
Moreover, if (i) and (ii) hold, then $\pi(G)$ is symmetric about $0$.
\end{lemma}

By virtue of Lemma \ref{-1 is in pi(G)}, it suffices to characterize those bipartite graphs $G$ such that $[0,1]\subset\pi(G)$.
For this purpose, we invoke the following result due to Schoenberg.
Observe that the $Q$-matrix makes sense for any metric space $\Gamma=(X,d)$ and hence so does the set $\pi(\Gamma)$.

\begin{theorem}[{cf.~\cite[Section 9.1]{PR2016B}}]\label{quadratic embedding}
For a metric space $\Gamma=(X,d)$, the following are equivalent:
\begin{enumerate}
\item $[0,1]\subset\pi(\Gamma)$.
\item $\Gamma$ has a quadratic embedding into a Hilbert space $\mathcal{H}$; i.e., there is a map $\theta:X\rightarrow\mathcal{H}$ such that $\|\theta(x)-\theta(y)\|^2=d(x,y)$ $(x,y\in X)$.
\item The function $d$ is a negative definite kernel on $X$; i.e.,
\begin{equation*}
	\sum_{x,y\in X} \!\overline{f(x)}f(y)\,d(x,y) \leqslant 0
\end{equation*}
for every $f\in C_0(X)$ such that $\sum_{x\in X}f(x)=0$.
\end{enumerate}
\end{theorem}

For two adjacent vertices $x$ and $y$, let
\begin{equation*}
	G(x,y)=\{z\in V:\partial(z,x)=\partial(z,y)-1\}.
\end{equation*}
In Section \ref{sec: proofs}, we next prove the following result, which gives a combinatorial characterization of graphs $G$ for which $\pi(G)=[-1,1]$:

\begin{theorem}\label{main theorem}
For a connected simple graph $G$ with at least two vertices, the following are equivalent:
\begin{enumerate}
\item $\pi(G)=[-1,1]$.
\item $G$ is isometrically embeddable into a hypercube (in a wide sense; see below).
\item $G$ is bipartite and has no quintuple $(x_1,x_2,x_3,x_4,x_5)$ of vertices such that $x_1\sim x_2$, $x_3\sim x_4$, $x_3\in G(x_1,x_2)$, $x_4\in G(x_2,x_1)$, and $x_5\in G(x_1,x_2)\cap G(x_4,x_3)$, where $\sim$ means adjacent.
\item $G$ is bipartite, and $G(x,y)$ is convex for any adjacent vertices $x$ and $y$; i.e., every geodesic between two vertices of $G(x,y)$ is entirely contained in $G(x,y)$.
\end{enumerate}
\end{theorem}

\noindent
The equivalence of (ii) and (iv) above is due to Djokovi\'{c} \cite{Djokovic1973JCTB}.
Here, we define the \emph{hypercube} $\mathcal{Q}_{\Sigma}$ on a nonempty set $\Sigma$ to be the graph having as vertices the finite subsets of $\Sigma$, where two vertices are adjacent if and only if one of them is obtained from the other by removing one element.
The distance between two vertices of $\mathcal{Q}_{\Sigma}$ agrees with the size of their symmetric difference.
Note that $\mathcal{Q}_{\Sigma}$ has a quadratic embedding into the Hilbert space $\ell^2(\Sigma)$.
An example of a quintuple in Theorem \ref{main theorem}\,(iii) is an induced $K_{2,3}$:
\begin{center}
\includegraphics[width=2.5cm,bb=0 13 300 183,clip]{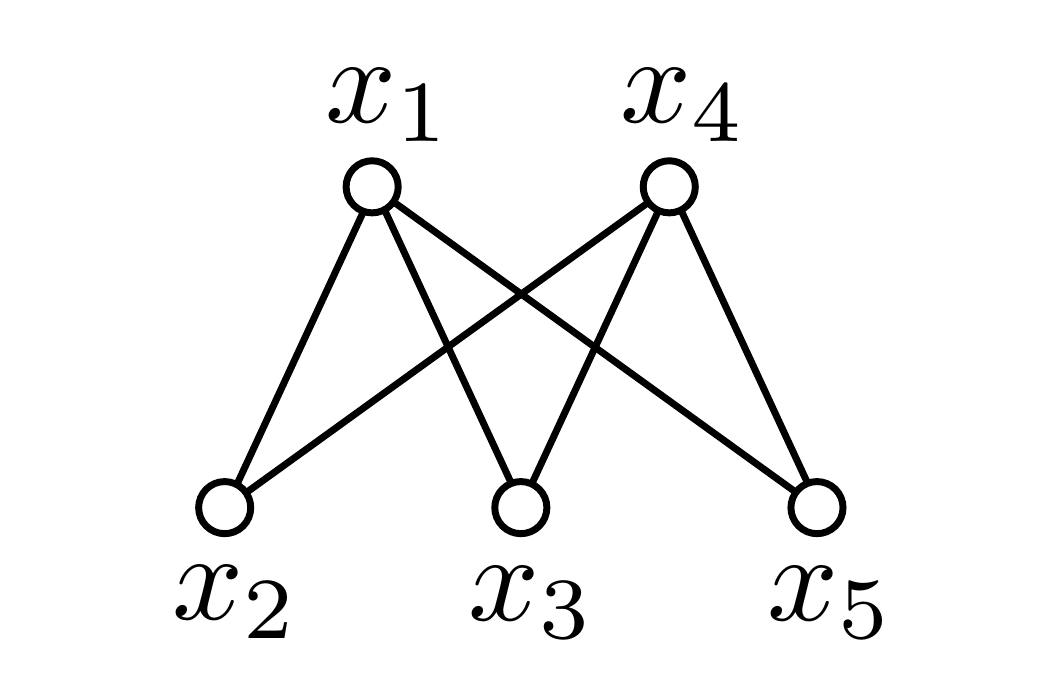}
\end{center}
From Theorem \ref{main theorem} it is immediate to see, for example, that $\pi(G)=[-1,1]$ whenever $G$ is a tree, a known result proved by many authors by different methods; see, e.g., \cite{Bozejko1989SM}.
We will discuss more examples in Sections \ref{sec: distance-regular graphs} and \ref{sec: Coxeter groups}.

In passing, Obata and Zakiyyah \cite{OZ2018EJGTA} defined the \emph{quadratic embedding constant} of $G$ as follows:
\begin{align*}
	\MoveEqLeft[1] \operatorname{QEC}(G) \\
	&=\sup\!\left\{ \sum_{x,y\in V} \!\overline{f(x)}f(y)\,\partial(x,y) : f\in C_0(V),\,\sum_{x\in V}|f(x)|^2=1,\,\sum_{x\in V}f(x)=0\right\}\!.
\end{align*}
By Theorem \ref{quadratic embedding}, $G$ admits a quadratic embedding into a Hilbert space if and only if $\operatorname{QEC}(G)\leqslant 0$, and hence this constant provides a new qualitative characteristic of graphs.
See \cite{Obata2022pre2} and the references therein for results on $\operatorname{QEC}(G)$.

\section{Proofs}\label{sec: proofs}

\begin{proof}[Proof of Lemma \ref{-1 is in pi(G)}]
(i)$\Rightarrow$(ii):
We begin by claiming that $G$ does not contain a triple $(x_1,x_2,x_3)$ of vertices such that $\partial(x_1,x_2)=\partial(x_1,x_3)$ and $x_2\sim x_3$.
Suppose, on the contrary, that such a triple exists, and write $i=\partial(x_1,x_2)=\partial(x_1,x_3)$.
Then, the principal submatrix of $Q_{-1}$ indexed by this triple is given by
\begin{equation*}
	\begin{bmatrix} 1 & (-1)^i & (-1)^i \\ (-1)^i & 1 & -1 \\ (-1)^i & -1 & 1 \end{bmatrix},
\end{equation*}	
which must be positive semidefinite.
However, this matrix has eigenvalues $2,2,-1$, so this is impossible.
The claim is proved.
We now show by induction that $G$ does not contain odd cycles.
$G$ does not contain $3$-cycles by the above claim.
Suppose that $x_0\sim x_1\sim\dots\sim x_{2t}\sim x_0$ is a $(2t+1)$-cycle in $G$, where $t\geqslant 2$.
We show by subinduction that $\partial(x_0,x_i)=i$ $(i=1,2,\dots,t)$.
Suppose that $\partial(x_0,x_i)<i$.
Then, since $\partial(x_0,x_i)\ne \partial(x_0,x_{i-1})=i-1$ by the above claim, we have $\partial(x_0,x_i)=i-2$ by the triangle inequality.
Let $y_1,y_2,\dots,y_{i-3}$ be vertices such that $x_0\sim y_1\sim y_2\sim\dots\sim y_{i-3}\sim x_i$.
Then by replacing $x_0\sim x_1\sim\dots\sim x_i$ by this sequence, we obtain a $(2t-1)$-cycle, which is a contradiction.
It follows that $\partial(x_0,x_i)=i$.
By the same reasoning, we also have $\partial(x_0,x_{t+1})=t+1$, but this is impossible since $\partial(x_0,x_{t+1})\leqslant t$ by $x_0\sim x_{2t}\sim\cdots\sim x_{t+2}\sim x_{t+1}$.
It follows that $G$ cannot contain a $(2t+1)$-cycle.
We have shown that $G$ does not contain odd cycles, so $G$ is bipartite.

(ii)$\Rightarrow$(i):
We first prove the last statement.
Fix a vertex $x_0$, and let $\Lambda$ be the diagonal matrix with $(x,x)$-entry
\begin{equation*}
	\Lambda_{x,x}=(-1)^{\partial(x_0,x)} \quad (x\in V).
\end{equation*}
Then we have $\Lambda Q_q\Lambda=Q_{-q}$.
To see this, note that
\begin{equation*}
	(\Lambda Q_q\Lambda)_{x,y}=(-1)^{\partial(x_0,x)+\partial(x_0,y)}q^{\partial(x,y)} \quad (x,y\in V),
\end{equation*}
and observe that $\partial(x_0,x)+\partial(x_0,y)+\partial(x,y)$ is even since $G$ is bipartite.
It follows that $Q_q\succcurlyeq 0$ if and only if $Q_{-q}\succcurlyeq 0$, and hence $\pi(G)$ is symmetric about $0$.
In particular, since $1\in \pi(G)$, we have $-1\in \pi(G)$.
\end{proof}

\begin{proof}[Proof of Theorem \ref{main theorem}]
(iii)$\Leftrightarrow$(iv):
Observe that a subset $U$ of $V$ is convex if and only if $z\in U$ for every vertex $z$ such that $x\sim z$ and $y\in G(z,x)$ for some $x,y\in U$.

(ii)$\Rightarrow$(i):
If $G$ is isometrically embeddable into the hypercube $\mathcal{Q}_{\Sigma}$, then $G$ has a quadratic embedding into $\ell^2(\Sigma)$.
By Lemma \ref{-1 is in pi(G)} and Theorem \ref{quadratic embedding}, we have (i).

(i)$\Rightarrow$(iii):
By Lemma \ref{-1 is in pi(G)}, $G$ is bipartite.
Suppose, on the contrary, that $G$ has a quintuple $(x_1,x_2,x_3,x_4,x_5)$ as in (iii), and write $i=\partial(x_1,x_3)$, $j=\partial(x_1,x_5)$, and $h=\partial(x_4,x_5)$.
Then, $\partial(x_2,x_3)=i+1$, $\partial(x_2,x_5)=j+1$, and $\partial(x_3,x_5)=h+1$.
Moreover, note that $\partial(x_1,x_4)=i\pm 1$ as $x_3\sim x_4$, but $i-1$ is ruled out since this would mean that $x_4$ lies on a geodesic between $x_1$ and $x_3$, implying $x_4\in G(x_1,x_2)$, a contradiction.
Finally, $\partial(x_2,x_4)=i$.
For five scalars $\xi_1,\xi_2,\xi_3,\xi_4,\xi_5\in\mathbb{R}$, we have
\begin{align*}
	\sum_{a,b=1}^5 \!\xi_a\xi_b\,\partial(x_a,x_b) = {} & 2\big( \xi_1\xi_2 +\xi_3\xi_4 + i(\xi_1\xi_3+\xi_2\xi_4) +(i+1)(\xi_1\xi_4+\xi_2\xi_3) \\[-3mm] & \quad + j\xi_1\xi_5 +(j+1)\xi_2\xi_5 +(h+1)\xi_3\xi_5 +h\xi_4\xi_5 \big).
\end{align*}
Set
\begin{equation*}
	\xi_1=\xi_4=-(j+h+2), \quad \xi_2=\xi_3=j+h, \quad \xi_5=4.
\end{equation*}
Then they sum to $0$, and
\begin{equation*}
	\sum_{a,b=1}^5 \!\xi_a\xi_b\,\partial(x_a,x_b)=8(i+1)>0.
\end{equation*}
But this is absurd since the distance $\partial$ is a negative definite kernel by Theorem \ref{quadratic embedding}.
Hence, (iii) follows.

(ii)$\Leftrightarrow$(iv):
See \cite{Djokovic1973JCTB}.
\end{proof}

\section{Distance-regular graphs}\label{sec: distance-regular graphs}

Distance-regular graphs are particularly well suited for the method of \emph{quantum decomposition} of the adjacency matrix, developed by Hora, Obata, and others; cf.~\cite{HO2007B}.
See also \cite{KOT2021SIGMA}.
Here, we generalize the definition of these graphs slightly, in the spirit of \cite{Zieschang1996B}.
Let us say that a connected simple graph $G=(V,E)$ is \emph{distance-regular} if there exist cardinal numbers $p_{i,j}^h$ $(i,j,h=0,1,2,\dots)$ such that, for any two vertices $x$ and $y$ with $\partial(x,y)=h$, we have
\begin{equation*}
	|\{z\in V:\partial(x,z)=i,\,\partial(z,y)=j\}|=p_{i,j}^h.
\end{equation*}
See \cite{BCN1989B,DKT2016EJC} for more information and recent updates on finite distance-regular graphs.
The finite distance-regular graphs that are isometrically embeddable into hypercubes were determined by Weichsel \cite{Weichsel1992DM}.
See also \cite{KS1994EJC}.
It follows that

\begin{theorem}
The finite distance-regular graphs $G$ for which $\pi(G)=[-1,1]$ are the finite hypercubes, the doubled Odd graphs, and the even cycles.
\end{theorem}

\noindent
For every positive integer $m$, the \emph{doubled Odd graph} $\tilde{O}_{m+1}$ can be defined as the subgraph of $\mathcal{Q}_{\Sigma}$ with $|\Sigma|=2m+1$ induced on the set of vertices which are $m$-subsets or $(m+1)$-subsets of $\Sigma$; cf.~\cite[Section 9.1D]{BCN1989B}.

Examples of infinite distance-regular graphs that are isometrically embeddable into hypercubes are the infinite hypercubes and the homogeneous trees.
It would be a problem of interest to determine all such graphs.

\section{Cayley graphs on Coxeter groups}\label{sec: Coxeter groups}

Cayley graphs on Coxeter groups were also studied in detail in the context of quantum probability theory; cf.~\cite{HO2007B}.
Recall that a \emph{Coxeter system} is a pair $(W,S)$ consisting of a group $W$ and a set of generators $S\subset W$, subject only to the relations
\begin{equation*}
	(ss')^{m(s,s')}=1 \quad (s,s'\in S),
\end{equation*}
where $m(s,s)=1$ and $m(s,s')=m(s',s)\in\{2,3,\dots\}\cup\{\infty\}$ if $s\ne s'$.
When $m(s,s')=\infty$, we understand that $ss'$ is of infinite order.
We use standard facts about Coxeter groups from \cite[Chapter 5]{Humphreys1990B}.
Let $G$ be the Cayley graph on $W$ with generating set $S$, so $G$ has vertex set $W$, and two vertices $x$ and $y$ are adjacent if and only if $yx^{-1}\in S$.
The $s\in S$ have order two, so $G$ is simple.
Clearly, $G$ is connected.
The length (of a reduced expression) of $x\in W$ is denoted by $\ell(x)$.
Note that $\partial(x,y)=\ell(yx^{-1})$ $(x,y\in W)$.
Let $\mathbb{V}$ be an $\mathbb{R}$-vector space with a basis $\{\alpha_s:s\in S\}$, equipped with a symmetric bilinear form $B$ given by
\begin{equation*}
	B(\alpha_s,\alpha_{s'})=-\cos\!\left(\frac{\pi}{m(s,s')}\right) \quad (s,s'\in S),
\end{equation*}
where the RHS is interpreted as $-1$ when $m(s,s')=\infty$.
By the \emph{geometric representation} $W\rightarrow \mathrm{GL}(\mathbb{V})$, the $s\in S$ act on $\mathbb{V}$ as
\begin{equation*}
	sv=v-2B(\alpha_s,v)\alpha_s \quad (v\in\mathbb{V}).
\end{equation*}
Let $\Phi$ be the set of \emph{roots} $x\alpha_s$ $(x\in W,\,s\in S)$.
Every root $\alpha\in \Phi$ is written as a linear combination of the $\alpha_s$, where either all coefficients are nonnegative or all coefficients are nonpositive.
The root $\alpha$ is \emph{positive} for the former and \emph{negative} for the latter.
Let $\Pi$ be the set of positive roots.
Then $\Phi=\Pi\sqcup(-\Pi)$, and we have
\begin{equation}\label{length}
	\ell(x) = |\Pi\cap x^{-1}(-\Pi)| \quad (x\in W).
\end{equation}
Now, let $\theta(x)=\Pi\cap x^{-1}(-\Pi)$ $(x\in W)$.
Then $|\theta(x)|=\ell(x)<\infty$, so the $\theta(x)$ are vertices of the hypercube $\mathcal{Q}_{\Pi}$.
Moreover, for any two vertices $x$ and $y$, we have
\begin{align*}
	\theta(x)\setminus\theta(y) &= \Pi\cap x^{-1}(-\Pi)\cap y^{-1}\Pi = -x^{-1}\big(x(-\Pi)\cap\Pi\cap xy^{-1}(-\Pi)\big), \\
	\theta(y)\setminus\theta(x) &= \Pi\cap y^{-1}(-\Pi)\cap x^{-1}\Pi = x^{-1} \big(x\Pi\cap\Pi\cap xy^{-1}(-\Pi)\big),
\end{align*}
so that
\begin{align*}
	|\theta(x)\triangle\theta(y)| &= |x(-\Pi)\cap\Pi\cap xy^{-1}(-\Pi)| + |x\Pi\cap\Pi\cap xy^{-1}(-\Pi)| \\
	&= |\Pi\cap xy^{-1}(-\Pi)|,
\end{align*}
and hence $|\theta(x)\triangle\theta(y)|=\ell(yx^{-1})=\partial(x,y)$.
It follows that the map $x\mapsto \theta(x)$ $(x\in W)$ is an isometric embedding of $G$ into $\mathcal{Q}_{\Pi}$.

\begin{theorem}
Let $(W,S)$ be a Coxeter system and let $G$ be the Cayley graph on $W$ with generating set $S$.
Then we have $\pi(G)=[-1,1]$.
\end{theorem}

\noindent
This result is not new.
Bo\.{z}ejko and Szwarc \cite{BS2003P} proved it using \eqref{length} and Theorem \ref{quadratic embedding}.
What is new here is an explicit mention of an isometric embedding into a hypercube.



\section*{Acknowledgments}
The author thanks Marek Bo\.{z}ejko for suggesting the problem discussed in this note.
Thanks are also due to the anonymous referee for offering valuable comments.
The author was supported by JSPS KAKENHI Grant Numbers JP19H01789 and JP20K03551.
This work was also partially supported by the Research Institute for Mathematical Sciences at Kyoto University.


\end{document}